\documentclass[12pt,a4paper,notitlepage]{article}

\usepackage{bbm}

\newcommand \NN{\mathbbm{N}}

\begin{document}

\title{Matheme and mathematics. On the main concepts of the philosophy of Alain Badiou.}
\author{Maciej Malicki \\
Department of Mathematics and Mathematical Economics, \\
Warsaw School of Economics, \\
mamalicki@gmail.com}
\date{Jan 19, 2014}

\maketitle

\begin{abstract}
In this paper, I present a critical discussion of mathematical arguments employed
in the philosophy of event of Alain Badiou. On the basis of ``Being and Event'' as
well as his other writings, I analyze the main notions of his philosophy such as the
indiscernible, the undecidable, and the unnameable. The focus of my analysis is
both on their mathematical consistency, and their philosophical consequences.
I argue that the mathematical approach developed by Badiou is seriously defective,
and, as a result, that it cannot serve as an ontological basis for the concept of event
as presented in ``Being and Event''.
\end{abstract}

\section{Introduction}

 ``Being and event'' is an ontological theory of the \emph{effect} of event, it is a
science of the possibility of novelty that breaks with being’s determinants.
Its existence can never be decided once and for all, because the occurrence
of an event is always an \emph{intervention}: its ultimate source lies in a decision
of the subject who is constituted in this act. As a result, a post-evental \emph{truth}
enters being. It is \emph{udecidable} and \emph{indiscernible}, which means that it cannot
be grounded nor described in terms of authoritative knowledge. It modifies
the situation by overruling some elements of knowledge, and making other
\emph{veridical}. Even though a truth cannot be accessed in its entirety, the subject
supporting a \emph{procedure of truth} can decide that something belongs to its realm.
These \emph{gestures of fidelity} give rise to the \emph{subject language}, which forces
statements describing the state of affairs after the occurrence of an event.

These are the core concepts of Badiou’s philosophy of event in a nutshell.
However --- and that is perhaps the most striking feature of ``Being and
Event'' --- all of its fundamental categories have a mathematical character,
and the main theses of this book derive from subtle and quite technical
considerations based on the area of mathematics called set theory. Eventually,
the science of event turns out to be the result of an investigation into
set theory, understood as ontology, or rather proved to be ontology, i.e. the
theory of ``being qua being''.

In this paper, I am mostly concerned with the consistency of mathematical
aspects of Badiou’s analysis, as well as some of its philosophical consequences
in certain strictly defined key aspects. It is an important reservation,
as a large part of Badiou’s thought can be presented without referring to any
technical terms. Many of its categories --- such as event itself, subject,
fidelity, the void or excess --- are well rooted in the language of contemporary
philosophy.

Taking a broader perspective, I am interested in the postulated deep
ontological structure behind four fundamental domains of truth, which are
science, politics, art and love. Of course, we do not need any special doctrine
in order to recognize the groundbreaking character of basic examples
of event discussed by Badiou --- Mallarme’s poetry, relativistic physics, the
communist revolution or dodecaphony, as well as some similarities between
them --- a break involved in their occurrence, the resulting new paradigm,
a peculiar indeterminateness of their status. On the other hand, some of
these examples do not seem to have that much in common with others, at
least on the surface of things. In order to acknowledge the relevance of
Badiou’s position, more is needed than a simple rephrasing of well known
ideas in a new language: what is required is a convincing argument showing
that the ontological theory of event is an autonomous interpretative tool,
independent of specific subject matters. Moreover, it is rather hard to imagine
an argument of this kind that would not contain a proof that the mathematical
part of Badiou’s thought is consistent and agrees with his fundamental
philosophical aims. The formalized doctrine presented by Badiou not only
allows for a realization of this postulate but it explicitly requires it.

Finally, one could describe this article as a contribution to the reflection
on the concept of matheme. In ``Being and Event'', a matheme is understood
as a philosophical idea subjected to rigors of deduction, and opposed to the
pre-platonic poem. However, the very term `matheme' comes from late
writings of Lacan, which are an important reference point for Badiou.
There it denotes mathematical objects — such as the Boromean knot or
the Klein bottle --- allowing to grasp the order of the real extending beyond
the reach of language. According to this understanding, a matheme would
be a place in philosophy, where mathematics attains an autonomous status,
completing philosophical discourse, and generating statements that are
binding for it.

My assessment of the method employed in ``Being and Event'' is negative.
After providing the first part of Badiou’s doctrine, that is the identification
of ontology with set theory, and outlining mathematical foundations of the
theory of event, I present a critical discussion of four key concepts of this
theory: the indiscernible, the undecidable, the unnameable and the evental
site. The conclusion is that their mathematical structure and its implications
for Badiou’s philosophy turn out not to meet expectations. They are ill-designed,
and this leads to mathematical inconsistencies as well as philosophical
consequences that contradict Badiou’s intentions. Far from rejecting the
philosophical substance of these categories, I claim that their mathematical
side may serve as an inspiring metaphor or analogy, but it has clearly
defined bounds of meaningful interpretation.

\section{Mathematics as ontology}
The edifice of Badiou’s philosophy of ``being qua being'' is founded on a
thesis establishing the identity of mathematics and ontology. In order to
accomplish this task, Badiou, in the very first words of ``Being and Event''
goes back to the old Parmenidean dispute on the status of the one and the
multiple. He decides it by giving priority to the multiple, while considering
the one as an effect of operation of imposing structure upon being, organizing
the ‘primordial’ \emph{inconsistent multiplicity}. Inconsistent means here
exactly this: ``without any unifying principle''.

Being such as it presents itself --- and ``no access to being is offered to
us except presentations'' [BE, p. 27] --- is a \emph{consistent multiplicity}. It is a
multiplicity because ``if the one is not reciprocal with being, the multiple,
however, is reciprocal with presentation.'' [BE, s. 28]; also, it is consistent
because a presented being is always the result of certain organizing principle,
certain \emph{operation of counting-as-one}. The inconsistent multiplicity is a
\emph{subtractive} basis of an already structured presentation --- it is possible to
be discerned only retroactively, as that ‘something’ on which a count-as-one
operated [BE, p. 25]. From the perspective of presentation it is only the void
because the effect of structure encompasses everything without exception. In
this manner, the inconsistent multiplicity ``sutures presentation to being''
[BE, p. 55]; the void is the name of being [BE, p. 56]. This establishment
of a consistent multiplicity resulting from some operation of counting-as-one,
a combination of the earlier inconsistent and later consistent multiplicity,
is called a \emph{situation}.

It is the inconsistent multiplicity that forms the domain of ontology.
``Ontology can be solely the theory of inconsistent multiplicities as such.
‘As such’ means that what is presented in the ontological situation is the
multiple without any other predicate than its multiplicity.'' [BE, p. 28] How
is the presentation of being regarded as multiplicity --- the presentation of
presentation --- possible if ``being has no structure'' [BE, p. 27]? What
requirements does such a situation need to satisfy? Firstly, a multiple without
one is a multiple that consists only of multiplicities. Secondly, the operation
of counting-as-one of the ontological situation cannot be anything more
than a collection of conditions through which the multiple can be recognized
as multiple. Finally, ontology must be a theory of ``the suture of presentation
to being'', that is, a theory of the subtractive void of presentation. In other
words, it should derive the existence of its multiplicities only out of the void.

The science satisfying these postulates is set theory --- a mathematical
theory of the relation of belonging , together with appropriate axioms:
the axiom of extensionality, regularity, pairing, union, infinity, power set,
and the axiom schema of comprehension and replacement. First of all, the
universe of set theory contains only sets, so set theory presents only multiplicities
consisting of multiplicities. Moreover, set theory is an axiomatic
theory, and axioms are formal rules that do not explicitly define objects they
refer to. Badiou says: ``an axiomatic presentation consists, on the basis of
non-defined terms, in prescribing the rule for their manipulation.'' [BE, p. 29]
Sets \emph{de facto} satisfy axioms of set theory, however the axioms themselves
do not form the criteria of being a set. In other words, what axioms determine
about sets is that they are multiplicities, nothing more. Finally, the
unpresented subtractive void of situation appears in set theory in the form of
the empty set, which is a basic building block of all sets. It must be a set
because every element presented in the ontological situation is the effect of
its operation of counting-as-one, however in fact, it is a ``‘multiple’ which is
neither one nor multiple, being the multiple of nothing, and therefore, as far
as it is concerned, presenting nothing in the form of the multiple, no more than
in the form of the one. This way ontology states that presentation is certainly
multiple, but that the being of presentation — that which is presented —
being void, is subtracted from the one/multiple dialectic.'' [BE, p. 59]

\section{Belonging, inclusion, and the impasse of being}

Let us take a closer look at the mathematical structure of key categories of
``Being and Event''. Ontologically, situations --- as it has been already said
--- take the form of sets: multiplicities whose one is nothing more than the
unity of the elements they consist of. The two basic types of relations
between sets are \emph{belonging}, that is being an element of a set, being \emph{presented
in a situation}, and \emph{inclusion}, that is, being a subset of a set, or being
\emph{represented in a situation}. Even though the latter relation is formally reducible
to the former one --- $y$ is a subset of $x$ if and only if every element of
$y$ is an element of $x$ --- their properties differ considerably. The tension
between belonging and inclusion is fundamental in so far as this is where
the ``impasse of being'' arises, opening up a situation to the interventional
occurrence of an event. More light can be shed on this tension with the help
of the concept of power set.

The power set $p(x)$ of $x$ --- the \emph{state of situation} $x$ or its \emph{metastructure}
[BE, p. 94] --- is defined as the set of all subsets of $x$. Now, basic relations
between sets can be expressed as the following relations between sets and
their power sets. If for some $x$, every element of $x$ is also a subset of $x$,
then $x$ is a subset of $p(x)$, and $x$ can be reduced to its power set. Conversely,
if every subset of $x$ is an element of $x$, then $p(x)$ is a subset of $x$, and the
power set $p(x)$ can be reduced to $x$. Sets that satisfy the first condition are
called \emph{transitive}. For obvious reasons the empty set is transitive; other
examples of transitive sets can also be easily found. However, the second
relation never holds. The mathematician Georg Cantor proved that not
only $p(x)$ can never be a subset of $x$, but in some fundamental sense it is
strictly larger than $x$. On the other hand, axioms of set theory do not determine
the extent of this difference. Badiou says that it is an ``excess of
being'', an excess that at the same time is its impasse.

In order to explain the mathematical sense of this statement, one needs to
recall the notion of \emph{cardinality}, which clarifies and generalizes the common
understanding of quantity. We say that two sets $x$ and $y$ have the same
cardinality if there exists a function defining a one-to-one correspondence
between elements of $x$ and elements of $y$. For finite sets, this definition agrees
with common intuitions: if a finite set $y$ has more elements than a finite set
$x$ (say, $y$ has $10$ elements and $x$ has $7$ elements), then regardless of how elements
of $x$ are assigned to elements of $y$, something (some $3$ elements) will
be left over in $y$ — precisely because it is larger. In particular, if $y$ contains
$x$ \emph{and} some other elements, then $y$ does not have the same cardinality as $x$.
This seemingly trivial fact is not always true outside of the domain of finite
sets. To give a simple example, the set of all natural numbers $\NN$ contains
quadratic numbers, that is, numbers of the form $n^2$, as well as some other
numbers but the set of all natural numbers, and the set of quadratic numbers
have the same cardinality. The correspondence witnessing this fact assigns
to every number $n$ a unique quadratic number, namely $n^2$.

Counting finite sets has always been done via natural numbers $0, 1, 2, \ldots$
In set theory, the concept of such a canonical measure can be extended to
infinite sets, using the notion of \emph{cardinal numbers}. Without getting into details
of their definition, let us say that the series of cardinal numbers begins with
natural numbers, which are directly followed by the number $\omega_0$, that is, the
size of the set of all natural numbers $\NN$, then by $\omega_1$, the first uncountable
cardinal numbers, etc. The hierarchy of cardinal numbers has the property
that every set $x$, finite or infinite, has cardinality (i.e. size) equal to exactly
one cardinal number $\kappa$. We say then that $\kappa$ is the cardinality of $x$.

It is not hard to prove that the cardinality of the power set $p(x)$ is $2^n$ for
every finite set $x$ of cardinality $n$. However, something quite paradoxical
happens when infinite sets are considered. Even though Cantor’s theorem
does state that the cardinality of $p(x)$ is always larger than $x$ --- similarly
as in the case of finite sets --- axioms of set theory \emph{never} determine the exact
cardinality of $p(x)$. Moreover, one can formally prove that there exists no
proof determining the cardinality of the power sets of any given infinite set.
There is a general method of building models of set theory, discovered by
the mathematician Paul Cohen, and called \emph{forcing}, that yields models, where
--- depending on construction details --- cardinalities of infinite power sets
can take different values. Consequently, quantity --- ``a fetish of objectivity''
[BE, p. 83] as Badiou calls it --- does not define a measure of being but it
leads to its impasse instead. It reveals an undetermined gap [BE, p. 83], where
an event can occur --- ``that-which-is-not being-qua-being'' [BE, p. 184].

\section{Forcing, truth and the place of the subject}

In order to make the exposition more accessible, let us consider only the
power set of the set $\NN$ of all natural numbers, which is the smallest infinite
set --- the countable infinity. Simplifying things slightly, the argument proceeds
as follows. By a model of set theory we understand a set in which
--- if we restrict ourselves to its elements only --- all axioms of set theory
are satisfied. It follows from G\"{o}del’s completeness theorem that as long as
set theory is consistent, no statement which is true in some model of set
theory can contradict logical consequences of its axioms. If the cardinality
of $p(\NN)$ was such a consequence, there would exist a cardinal number $\kappa$
such that the sentence ‘the cardinality of $p(\NN)$ is $\kappa$’ would be true in all the
models. However, for every cardinal $\kappa$ the technique of forcing allows for
finding a model $M$ where the cardinality of $p(\NN)$ is not equal to $\kappa$. Thus,
for no $\kappa$, the sentence ‘the cardinality of $p(\NN)$ is $\kappa$’ is a consequence of the
axioms of set theory, i.e. they do not decide the cardinality of $p(\NN)$.

The starting point of forcing is a model $M$ of set theory --- called the
\emph{ground model} --- which is countably infinite and transitive (this is a crucial
assumption). As a matter of fact, the existence of such a model cannot be
proved but it is known that there exists a countable and transitive model for
every \emph{finite} subset of axioms. As far as the logic of the construction is
concerned, in particular the decidability of sentences obtained by forcing,
this difference does not play any role.

A characteristic subtlety can be observed here. From the perspective of
‘an inhabitant of the universe’, that is, if all the sets are considered, the
model $M$ is only a small part of this universe. It is deficient in almost every
respect; for example all of its elements are countable, even though the existence
of uncountable sets is a consequence of the axioms of set theory.
However, from the point of view of an ‘inhabitant of $M$’, that is, if elements
outside of $M$ are disregarded, everything is in order. Some of $M$ because
in this model there are no functions establishing a one-to-one correspondence
between them and $\omega_0$. One could say that $M$ “simulates” the properties
of the whole universe.

The main objective of forcing is to build a new model $M[G]$ based on
$M$, which contains $M$, and satisfies certain additional properties. The model
$M[G]$ is called the \emph{generic extension} of $M$. In order to accomplish this goal,
a particular set is distinguished in $M$ — its elements are referred to as
conditions — which will be used to determine basic properties of the
generic extension. In case of the forcing that proves the undecidability of
the cardinality of $p(\NN)$, the set of conditions codes finite fragments of a function
witnessing the correspondence between $p(\NN)$ and a fixed cardinal $\kappa$.

In the next step, an appropriately chosen set $G$ is added to $M$ as well as
other sets that are indispensable in order for $M[G]$ to satisfy the axioms of
set theory. This set --- called \emph{generic} --- is a subset of the set of conditions
that always lays outside of $M$. The construction of $M[G]$ is exceptional in
the sense that its key properties can be described and proved using $M$
only, or just the conditions, thus, without referring to the generic set. This
is possible for three reasons. First of all, every element $x$ of $M[G]$ has a
“name” existing already in $M$ (that is, an element in $M$ that codes $x$ in some
particular way). Secondly, based on these names, one can design a language
called the \emph{forcing language} or --- as Badiou terms it --- the \emph{subject language}
that is powerful enough to express every sentence of set theory referring to
the generic extension. Finally, it turns out that the validity of sentences of
the forcing language in the extension $M[G]$ depends on the set of conditions:
the conditions \emph{force} validity of sentences of the forcing language in a precisely
specified sense. As it has already been said, the generic set $G$ consists
of some of the conditions, so even though $G$ is outside of $M$, its elements
are in $M$. Recognizing which of them will end up in $G$ is not possible for
an inhabitant of $M$, however in some cases the following can be proved:
provided that the condition $p$ is an element of $G$, the sentence $S$ is true in
the generic extension constructed using this generic set $G$. We say then that
$p$ forces $S$.

In this way, with an aid of the forcing language, one can prove that every
generic set of the Cohen forcing codes an entire function defining a one-to-one
correspondence between elements of $p(\NN)$ and a fixed (uncountable) cardinal number
--- it turns out that \emph{all} the conditions force the sentence stating this property
of $G$, so regardless of which conditions end up in the generic set, it is
always true in the generic extension. On the other hand, the existence of a
generic set in the model $M$ cannot follow from axioms of set theory, otherwise
they would decide the cardinality of $p(\NN)$.

The method of forcing is of fundamental importance for Badiou’s philosophy.
The event escapes ontology; it is ``that-which-is-not-being-qua-being'',
so it has no place in set theory or the forcing construction. However,
the post-evental truth that enters, and modifies the situation, is presented by
forcing in the form of a generic set leading to an extension of the ground
model. In other words, the situation, understood as the ground model $M$, is
transformed by a post-evental truth identified with a generic set $G$, and
becomes the generic model $M[G]$. Moreover, the knowledge of the situation
is interpreted as the language of set theory, serving to discern elements
of the situation; and as axioms of set theory, deciding validity of statements
about the situation. Knowledge, understood in this way, does not decide the
existence of a generic set in the situation nor can it point to its elements
(this property of truth will be thoroughly discussed later). A generic set is
always \emph{undecidable} and \emph{indiscernible}.

Therefore, from the perspective of knowledge, it is not possible to establish,
whether a situation is still the ground-model or it has undergone a
generic extension resulting from the occurrence of an event; only the subject
can interventionally decide this. And it is only the subject who decides
about the belonging of particular elements to the generic set (i.e. the truth).
A \emph{procedure of truth} or \emph{procedure of fidelity} [BE, p. 329] supported in this
way gives rise to the subject language. It consists of sentences of set theory,
so in this respect it is a part of knowledge, although the \emph{veridicity} of the
subject language originates from decisions of the faithful subject. Consequently,
a procedure of fidelity forces statements about the situation as it is
after being extended, and modified by the operation of truth.

\section{Mathemes of the undecidable and of the evental site}

According to Badiou, the undecidable truth is located beyond the boundaries
of authoritative claims of knowledge. At the same time, undecidability
indicates that truth has a post-evental character: ``the heart of the truth is
that the event in which it originates is undecidable'' [BE, p. 221]. Badiou
explains that, in terms of forcing, undecidability means that the conditions
belonging to the generic set force sentences that are not consequences of
axioms of set theory. However, one also needs to answer the question about
the role played by axioms in the structure of \emph{historical} situations. If in the
domains of specific languages (of politics, science, art or love) the effects
of event are not visible, the content of ``Being and Event'' is an empty
exercise in abstraction: even science --- perhaps excluding some entirely
formalized areas of theoretical physics --- let alone art or love --- cannot
for obvious reasons be exhaustively described solely in terms of the relation
of belonging. Anyway, it is doubtful that --- to consider just one example
--- the status of the French revolution is different from the status of the
absolute monarchy preceding it as far as the axioms of set theory are concerned.
In other words, most likely either all historical facts are decidable or
none of them is. Both possibilities lead to a trivial notion of the undecidable.

Judging by numerous examples discussed by Badiou, it seems that he
distances himself form such a narrow interpretation of the function played
by axioms. He rather regards them as collections of basic convictions
that organize situations, the conceptual or ideological framework of a historical
situation. For example, the nature of politics in Rousseau’s writings
is formulated in the following way: ``The major axiom is that in order to
definitely have the expression of the general will, [there must] be no partial
society in the State'' [BE, p. 348]. This approach is also indicated in the
only part of ``Being and Event'' which considers that issue in general terms:
``Let us agree that a proposition is singular (\ldots) if, within a historically
structured mathematical situation, it implies many other significant propositions,
yet it cannot itself be deduced from the axioms which organize the
situation. (\ldots) Say that A is this proposition. (\ldots) An event, named by an
intervention, is then, at the theoretical site indexed by the proposition $A$, a
new apparatus, demonstrative or axiomatic, such that A is henceforth clearly
admissible as a proposition of the situation.'' [BE, p. 246] Accordingly, the
undecidability of a truth would consist in transcending the theoretical
framework of a historical situation or even breaking with it in the sense that
the faithful subject accepts beliefs that are impossible to reconcile with the
old mode of thinking.

A clear illustration of the effect of event which in Badiou’s opinion
results in breaking with the determinants of the old paradigm, rather than
just in moving beyond them, is the birth of relativistic physics: ``After
Einstein’s texts of 1905, if I am faithful to their radical novelty, I cannot
continue to practice physics within its classical framework'' [Ethics, p. 42].
The novelty of relativistic physics cannot be reduced to a mere substitution
of certain equations with different, more precise ones, because it gives rise
to a completely new understanding of fundamental physical categories such
as space, time, reference point or motion. For profound reasons classical
mechanics rules out --- instead of simply not deciding it --- the very possibility
of the theory of relativity emerging within its own conceptual framework.
Similarly, the French revolution and communism --- essential examples of
the effect of event --- violently rupture the historical order preceding them.
The French revolution overthrew the king and established the sovereignty
of the people. And if something can be said with certainty about any real
or imaginary realization of the communist idea, it is that it definitely abolishes
the ``capital-parliamentarism''.

However, if one consequently identifies the effect of event with the structure
of the generic extension, they need to conclude that these historical
situations are by no means the effects of event. This is because a crucial
property of every generic extension is that axioms of set theory remain
valid within it. It is the very core of the method of forcing, stated in the
Theorem of the Generic Model [Jech, Th. 14.5]. Without this assumption,
Cohen’s original construction would have no \emph{raison d’etre} because it would
not establish the undecidability of the cardinality of infinite power sets. Let
us say this once more: every generic extension satisfies axioms of set theory.
In reference to historical situations, it must be conceded that a procedure
of fidelity may modify a situation by forcing undecidable sentences, nonetheless
it never overrules its organizing principles.

From the point of view of the generic theory of truth, some hypothetical
type of social democracy might be considered as the effect of event. It
would abolish chaos and inequalities, resulting from mechanisms of democratically
controlled market economy, by the operation of a new idea transfiguring
the nature of these mechanisms from within. As a religious event,
transgressing the Law without literally breaching it, one could probably
point to Messianic Judaism or Protestantism. Another interesting case is a
theory of the Danish astronomer Tycho de Brahe who in the 16th century
proposed a solution that allowed for keeping the empirical advantages of
the heliocentric model, while letting the Earth stay in the center of the
Universe. In this conception, all the planets revolve around the Sun, except
for Earth, which is encircled by the Sun. In terms of kinetics, that is, if the
force of gravity --- unknown at that time --- is disregarded, de Brahe’s model
is entirely equivalent to the Copernican one.

Another notion which cannot be located within the generic theory of
truth without extreme consequences is \emph{evental site}. An evental site --- an
element ``on the edge of the void'' [BE, p. 175] --- opens up a situation to
the possibility of an event [BE, p. 179]. Ontologically, it is defined as ``a
multiple such that none of its elements are presented in the situation'' [BE,
p. 175]. In other words, it is a set such that neither itself nor any of its
subsets are elements of the state of the situation. As the double meaning of
this word indicates, the state (\emph{\'{e}tat}) in the context of historical situations
takes the shape of the State (\emph{\'{E}tat}) [BE, p. 104]. A paradigmatic example of
a historical evental site is the proletariat --- ``entirely dispossessed, and
absent from the political stage'' [Ethics, p. 69].

The existence of an evental site in a situation is a necessary requirement
for an event to occur. Badiou is very strict about this point: ``we
shall posit once and for all that there are no natural events, nor are there
neutral events'' [BE, p. 178] --- and it should be clarified that situations
are divided into natural, neutral, and those that contain an evental site.
The very matheme of event --- its formal definition is of no importance
here — is based on the evental site [Ethics, p. 179]. The event raises the
evental site to the surface, making it represented on the level of the state of
the situation. Moreover, a novelty that has the structure of the generic set
but it does not emerge from the void of an evental site, leads to a \emph{simulacrum
of truth} [Ethics, p. 72], which is one of the figures of Evil [Ethics,
p. 87]. An example of utterly destructive effects of a simulacrum of truth is
the Nazi revolution whose source was the ``plenitude'' of the German people
[Ethics, p. 73].

However, if one takes the mathematical framework of Badiou’s concept
of event seriously, it turns out that there is no place for the evental site
there --- it is forbidden by the assumption of transitivity of the ground
model $M$. This ingredient plays a fundamental role in forcing, and its
removal would ruin the whole construction of the generic extension. As it
has already been mentioned, transitivity means that if a set belongs to $M$,
all its elements also belong to $M$. However, an evental site is a set none of
whose elements belongs to $M$. Therefore, contrary to Badiou’s intentions,
there cannot exist evental sites in the ground model. Using Badiou’s
terminology one can say that forcing may only be the theory of the simulacrum
of truth.

\section{The mathemes of the indiscernible and unnameable}

``Thought is nothing other than the desire to finish with the exorbitant
excess of the state'' [BE, p. 282]. Since Cantor’s theorem implies that this
excess cannot be removed or reduced to the situation itself, the only way
left is to take control of it. A basic, paradigmatic strategy for achieving this
goal is to subject the excess to the power of language. Its essence has been
expressed by Leibniz in the form of the principle of indiscernibles: there
cannot exist two things whose difference cannot be marked by a describable
property [BE, p. 283]. In this manner, language assumes the role of a ``law
of being'' [BE, p. 283], postulating identity, where it cannot find a difference.
Meanwhile --- according to Badiou --- the generic truth is indiscernible:
there is no property expressible in the language of set theory that characterizes
elements of the generic set. Truth is beyond the power of knowledge,
only the subject can support a procedure of fidelity by deciding what belongs
to a truth. This key thesis is established using purely formal means, so it
should be regarded as one of the peak moments of the mathematical method
employed by Badiou. In order to assess its grounding and possible limitations,
one needs to analyze the matheme of the indiscernible as closely as possible.

To the reader's surprise, Badiou composes the indiscernible out of as many as three different mathematical notions. First of all, he decides that it corresponds to the concept of the inconstructible \cite[p. 355]{BE}. Later, however, he writes that ``a set $\delta$ is discernible (\ldots) if there exists (\ldots) an explicit formula $\lambda(x)$ (\ldots) such that 'belong to $\delta$` and 'have the property expressed by $\lambda(x)$` coincide'' \cite[p. 367]{BE}. Finally, at the outset of the argument designed to demonstrate the indiscernibility of truth he brings in yet another definition: ``let us suppose the contrary: the discernibility of $G$. A formula thus exists $\lambda(x,a_1,\ldots, a_n)$ with parameters $a_1 \ldots,a_n$ belonging to $M[G]$ such that for an inhabitant of $M[G]$ it defines the multiple $G$'' \cite[p. 386]{BE}. In short, discernibility is understood as:

\begin{enumerate}
\item constructibility
\item  definability by a formula $F(y)$ with one free variable and no parameters. In this approach, a set $a$ is definable if there exists a formula $F(y)$ such that $b$ is an element of $a$ if and only if $F(b)$ holds.
\item  definability by a formula $F(y,z_1 \ldots,z_n)$ with parameters. This time, a set $a$ is definable if there exists a formula $F(y,z_1, \ldots,z_n)$ and sets $a_1,\ldots,a_n$ such that after substituting $z_1=a_1,\ldots,z_n=a_n$, an element $b$ belongs to $a$ if and only if $F(b,a_1,\ldots,a_n)$ holds.
\end{enumerate}

Even though in ``Being and Event'' Badiou does not explain the reasons for
this variation, it clearly follows from his other writings (such as [Conditions,
p. 135]) that he is convinced that these notions are equivalent.
It should be emphasized then that this is not true: a set may be discernible
in one sense, but indiscernible in another. First of all, the last definition has
been included probably by mistake because it is trivial. Every set in $M[G]$
is discernible in this sense because for every set a the formula $F(y, x)$
defined as ‘$y$ belongs to $x$’ defines a after substituting $x = a$. Accepting this
version of indiscernibility would lead to the conclusion that truth is always
discernible, while Badiou claims that it is not so. In particular, the proof of
the indiscernibility of truth presented by Badiou on page 386 of ``Being and
Event'', based on this definition of indiscernibility, is incorrect\footnote{There are several flaws in the proof but the most important one is that Badiou wrongly assumes that an element defined by a formula with parameters in the ground model must be an element of the ground model as well. This is not true as the following discussion will show.}.

Is it not possible to choose the second option and identify discernibility
with definability by a formula with no parameters? After all, this notion is
most similar to the original idea of Leibniz — intuitively, the formula $F(y)$
expresses a property characterizing elements of the set defined by it. Unfortunately,
this solution does not warrant indiscernibility of the generic set
either. There are examples of generic extensions $M[G]$, where the generic
set is definable, or even definable in $M[G]$, by a formula with no parameters
(see [Enayat], and [Fuchs].) Therefore, in this approach truth may be seized
by knowledge replacing the intervening subject.

As a matter of fact, assuming that in ontology, that is, in set theory, discernibility
corresponds to constructibility, Badiou is right that the generic
set is necessarily indiscernible. However, constructibility is a highly technical
notion, and its philosophical interpretation seems very problematic. Let us
take a closer look at it.

The class of constructible sets --- usually denoted by the letter $L$ ---
forms a hierarchy indexed or ‘numbered’ by ordinal numbers. Without
getting into details of the definition of the ordinal number --- closely related
to that of cardinal number --- the inductive procedure of constructing
the constructible hierarchy goes as follows. The lowest level $L_0$ is simply
the empty set. Assuming that some level --- let us denote it by $L_\alpha$ --- has
already been constructed, the next level $L_{\alpha+1}$ is constructed by choosing all
subsets of L that can be defined by a formula (possibly with parameters)
\emph{bounded} to the lower level $L_\alpha$.

Bounding a formula to $L_\alpha$ means that its parameters must belong to $L_\alpha$ and that its quantifiers are restricted to elements of $L_\alpha$. For instance, the formula 'there exists $z$ such that $z$ is in $y$` simply says that $y$ is not empty. After bounding it to $L_\alpha$ this formula takes the form 'there exists $z$ in $L_\alpha$ such that $z$ is in $y$`, so it says that $y$ is not empty, and some element from $L_\alpha$ witnesses it. Accordingly, the set defined by it consists of precisely those sets in $L_\alpha$ that contain an element from $L_\alpha$.

After constructing an infinite sequence of levels (or --- strictly speaking --- a limit sequence of levels) the level directly above them all is simply the set of all elements constructed so far. For example, the first infinite level $L_\omega$ consists of all elements constructed on levels $L_0, L_1, L_2, \ldots$.

As a result of applying this inductive definition, on each level of the
hierarchy all the formulas are used, so that two distinct sets may be defined
by the same formula. On the other hand, only bounded formulas take part
in the construction. The definition of constructibility offers too little (because
only bounded formulas are accepted) and too much at the same time
(because many sets can be defined by one formula). This technical notion
resembles the Leibnizian discernibility only in so far as it refers to formulas.
In set theory there are more notions of this type though.

To realize difficulties involved in attempts to philosophically interpret
constructibility, one may consider a slight, purely technical, extension of
it. Let us also accept sets that can be defined by a formula $F(y, z_1, \ldots, z_n)$
with constructible parameters, that is, parameters coming from $L$. Such a
step does not lead further away from the common understanding of Leibniz’s
principle than constructibility itself: if parameters coming from lower
levels of the hierarchy are admissible when constructing a new set, why not
admit others as well, especially since this condition has no philosophical
justification?

Actually, one can accept parameters coming from an even more restricted class,  e.g., the class of ordinal numbers. Then we will obtain the notion of \emph{definability from ordinal numbers} $(OD)$. This minor modification of the concept of constructibility --- a  relaxation of the requirement that the procedure of construction has to be restricted to lower levels of the hierarchy --- results in drastic consequences. Kenneth McAloon (\cite{McAloon}) proved in 1972 that the generic set as well as all other elements of the generic extension may be definable from ordinal numbers. Therefore, replacing constructibility with a concept very similar to it ---  one essentially identical from the point of view of Badiou's philosophical motivations --- leads to an ontology allowing for the discernibility of truth.

Another example of such a heterogeneous composition is the matheme of \emph{unnameable}. In ``Being and Event'' this notion functions as a synonym of the indiscernible, however in Badiou's later writings ---  especially ``Conditions'' and ``Ethics'' --- it becomes a distinct, autonomous concept. Referring to Lacanian psychoanalysis, Badiou says: ``in the field determined by a situation and the generic becoming of its truth, a real is attested to by a term, a point, and only one, at which the power of truth is suspended.  There is only one term in relation to which no anticipatory hypothesis on the generic subset permits us to force a judgment, (\ldots) no naming is appropriate for this term. That is why I call it unnameable'' \cite[p. 141]{C}. Also, the unnameable understood in the same manner --- as a ``point that the truth cannot force'' \cite[p. 85]{E} --- witnesses the ``powerlessness of truth'' in Badiou's ethical thought. A totally powerful truth  ---  its subject language reaching all the elements of the situation ---  leads to a ``disaster'', which is another figure of Evil, along with the simulacrum of truth \cite[p. 85]{E}.  

Thus in set theory, the unnameable has a precisely defined form: it is the \emph{unique} element of the generic extension such that no condition forces a sentence of the forcing language referring to this element. Unfortunately, the mathematical content of this definition turns out to be empty ---  forcing unequivocally rules out the existence of the unnameable. This is because one of the basic features of the forcing language is (\cite[Th. 14. 7 ii) a)]{Jech}) that if no condition forces a given sentence $S$, then every condition forces the negation of $S$. In the generic extension, the unnameable does not exist. 

Yet Badiou's position is that mathematics indicates the possibility of the unnameable. He points to a construction due to Furkhen ``in which one term exists, and only one, which cannot receive a name in the sense that it cannot be identified by a formula in the language'' \cite[p. 143]{C}. A familiar substitution can be recognized here. Initially, the unnameable is defined in terms of the forcing language, however Furkhen's construction refers to definability by a formula. It should be emphasized that this object is not even a model of set theory \cite[p. 120]{C}; it requires a richer  language, and it does not satisfy the axioms of set theory. Therefore it cannot serve as a ground model or a generic extension --- forcing has simply no use in this framework. Moreover, the unnameable becomes the indiscernible again with an additional requirement of uniqueness. Is a truth a real of the situation then? Such a suggestion never appears in Badiou's writings. Otherwise, there always exist at least two indiscernible elements --- a truth and a real of the situation --- even though the unnameable is supposed to be unique. In presence of these overwhelming formal and interpretative difficulties, one should rather conclude that the matheme of unnameable is an ill-conceived concept.

\begin{center} 
***
\end{center}

The mathematical form of the undecidable implies that a truth never overrules
the fundamental conceptual framework of a situation from before the
occurrence of an event. This consequence remains in sharp contrast with
basic examples of the effect of event considered by Badiou. An evental site
cannot open up a situation to an evental truth because forcing rules out the
very existence of an evental site in ground models. As a result, the generic
set may possibly describe the structure of the simulacrum of truth --- using
Badiou’s terminology --- but not the structure of truth itself. Two two
remaining mathemes --- the indiscernible and the unnameable --- are built out
of mathematical concepts that do not fit together; moreover, their properties
are not as Badiou claims.

This list of problems can be extended. The ``name of the event'' is not
``prohibited by being'' [BE, p. 184] because the axiom of foundation may
be consistently replaced with its negation --- for example, the Aczel antifoundation
axiom is consistent with the remaining axioms of set theory.
Finite sets are not necessarily definable, even though Badiou tries to prove
the opposite [Conditions, p. 137], mistakenly identifying finite sets with
numerals. In the so called nonstandard models of set theory the standard
omega is finite and non definable. Generic sets are not coextensive with
unconstructible sets [Conditions, p. 135] --- for example, under the Ground
Axiom. It seems though that the discussion of the mathemes of indiscernible,
undecidable, unnameable, and evental site is sufficient to defend the
main thesis of this paper. Mathematics --- forced to accept compromises
going definitely too far --- responds with outcomes which are hostile to
fundamental philosophical motivations of Badiou’s doctrine. Despite some
points of convergence, his generic theory of truth and his philosophy of
event can coexist only at a price of selective and instrumental interpretation
of the mathematical component. Therefore one has to conclude that ``Being
and Event'' provides no grounding for a deep ontological structure behind
the realms of science, art, love and politics, and that the mathematical formulation
of the theory of event has no positive content.


\begin{thebibliography}{99}
\bibitem[BE]{BE} A. Badiou, \emph{Being and Event}, Continuum 2005.
\bibitem[Ethics]{E} A.Badiou, \emph{Ethics:An Essay on the Understanding of Evil}, Verso 2002.
\bibitem[Conditions]{C} A.Badiou, \emph{Conditions}, Continuum 2009.
\bibitem[Enayat]{Enayat}A.Enyat, Models of set theory with definable ordinals, \emph{Archive for Math. Logic}, vol. 44, 3 (2005), pp. 363-385.
\bibitem[Fuchs]{Fuchs} G.Fuchs, J.D.Hamkins, Degrees of rigidity for Souslin trees, \emph{J. Symb. Logic}, vol. 74, 2 (2009), pp. 423--454.
\bibitem[Jech]{Jech} T.Jech, \emph{Set theory}, Springer 2006.
\bibitem[McAloon]{McAloon} K.McAloon, Consistency results about ordinal defnability, \emph{Annals of Math. Logic}, vol. 2  (1971),  pp.449-467.
\end{thebibliography}
\end{document}